\documentclass[12pt,leqno]{article}
\usepackage{amsmath,amsfonts,fancybox,qsymbols}
\reversemarginpar

\begin{document}
%
%%%%%%%%%%%%%%%%%%%%%%%%%%%%%%%%%%%%%%%%%%%%%%%%%%%%%%%%%%%%%%%%%%%%%%%%
\newtheorem{theorem}{Theorem}
\newtheorem{proposition}[theorem]{Proposition}
\newtheorem{lemma}[theorem]{Lemma}
\newtheorem{assum}[theorem]{Assumption}
\newtheorem{corollary}[theorem]{Corollary}
\newtheorem{rem}[theorem]{Remark}
\newqsymbol{`P}{\mathbb{P}}
\newqsymbol{`E}{\mathbb{E}}
%%%%%%%%%%%%%%%%%%%%%%%%%%%%%
%%%%%%%%%%%%%%%%%%%%%%%%%%%%%
\def\d{\displaystyle}
\def\pn{\par\noindent}
\def\hold{${\cal C}^{\alpha}([0,T])$ }
\def\hsq{${\cal C}^{1/2}([0,T])$ }
\def\mod{mod\`ele }
\def\sar{{\cal S}^{\downarrow}}
%%%%%%%%%%%%%%%%%%%%%%%%%%%%%
\newdimen\AAdi%
\newbox\AAbo%
\font\AAFf=cmex10%   %ou cmex10
\def\AArm{\fam0 }%\tenrm}%
\def\AAk#1#2{\setbox\AAbo=\hbox{#2}\AAdi=\wd\AAbo\kern#1\AAdi{}}%
\def\AAr#1#2#3{\setbox\AAbo=\hbox{#2}\AAdi=\ht\AAbo\raise#1\AAdi\hbox{#3}
}%
%%%%%%%%%%%%%%%%%%%%%%%%%%%%%
%% alphabet blackboard
\def\BBa{{\AArm A}}%
\def\BBb{{\AArm I\!B}}%
\def\BBc{{\AArm C\AAk{-1.02}{C}\AAr{.9}{I}{\AAFf\char"3F}}}%
\def\BBd{{\AArm I\!D}}%
\def\BBe{{\AArm I\!E}}%
\def\BBf{{\AArm I\!F}}%
\def\BBg{{\AArm G\AAk{-1.02}{G}\AAr{.9}{I}{\AAFf\char"3F}}}%
\def\BBh{{\AArm I\!H}}%
\def\BBi{{\AArm I\!I}}%
\def\BBj{{\AArm J}}%
\def\BBk{{\AArm I\!K}}%
\def\BBl{{\AArm I\!L}}%
\def\BBm{{\AArm I\!M}}%
\def\BBn{{\AArm I\!N}}%
\def\N{{\AArm I\!N}}%
\def\BBo{{\AArm O\AAk{-1.02}{O}\AAr{.9}{I}{\AAFf\char"3F}}}%
\def\BBp{{\AArm I\!P}}%
\def\P{{\AArm I\!P}}%
\def \p{{\cal P}}
\def\BBq{{\AArm Q\AAk{-1.02}{Q}\AAr{.9}{I}{\AAFf\char"3F}}}%
\def\Q{{\AArm Q\AAk{-1.02}{Q}\AAr{.9}{I}{\AAFf\char"3F}}}%
\def\BBr{{\AArm I\!R}}%
\def\BBs{{\AArm S}}%
\def\BBt{{\AArm T\AAk{-.62}{T}T}}%
\def\BBu{{\AArm U\AAk{-1}{U}\AAr{.95}{I}{\AAFf\char"3F}}}%
\def\BBv{{\AArm V}}%
\def\BBw{{\AArm W}}%
\def\BBx{{\AArm X}}%
\def\BBy{{\AArm Y}}%
\def\BBz{{\AArm Z\!\!Z}}%
\def\BBone{{\AArm 1\AAk{-.8}{I}I}}%
\def\D{{\cal D}}
\def \\ { \cr }
\def\R{{\math R}}
\def\N{{\math N}}
\def\E{{\math E}}
\def\P{{\math P}}
\def\Z{{\math Z}}
\def\Q{{\math Q}}
\def\C{{\math C}}
\def \e{{\rm e}}
\def \f{{\cal F}}
\def \g{{\cal G}}
\def \L{{\cal L}}
\def \d{{\tt d}}
\def \k{{\tt k}}
\def \i{{\tt i}}
\def \p{{\cal P}}
\def \s{{\cal S}}
\newcommand{\ed}{\mbox{$ \ \stackrel{d}{=}$ }}
\def\ni{\noindent}
\def \be{\begin{eqnarray*}}
\def \ee{\end{eqnarray*}}
\def \ben{\begin{eqnarray}}
\def \een{\end{eqnarray}}
%%%%%%%%%%%%%%%%%%%%%%%%%%%%%%%%%%%%%%%%%%%%%%%%%%%%%%%%%%%%%%%%%%%%%%%
%%%%%%%%%%%%%%%%%%%%%%%%%%%%%%%%%%%%%%%%%%%%%%%%%%%%%%%%%%%%%%%%%%%%%%%%
\def\videbox{\mathbin{\vbox{\hrule\hbox{\vrule height1ex \kern.5em\vrule height1ex}\hrule}}}
\font\calcal=cmsy10 scaled\magstep1
\def\build#1_#2^#3{\mathrel{\mathop{\kern 0pt#1}\limits_{#2}^{#3}}}
\def\suPP#1#2{{\displaystyle\sup _{\scriptstyle #1\atop \scriptstyle #2}}}
\def\proDD#1#2{{\displaystyle\prod _{\scriptstyle #1\atop \scriptstyle #2}}}
\def \s{{\cal S}}
\def\R{{\math R}}
\def\proof{\noindent{\bf Proof:}\hskip10pt}
\def\QED{\hfill\vrule height 1.5ex width 1.4ex depth -.1ex \vskip20pt}
%%%%%%%%%%%%%%%%%%%%%%%%%%%%%%%%%%%%%%%%%%%%%%%%%%%%%%%%%%%%%%%%%%%%%%%%
\setcounter{page}{1}
\setcounter{section}{0}
\title{Asymptotical behaviour of the presence probability\\ in branching random walks and fragmentations}
\maketitle
\vskip 0.5cm
\centerline{\Large \bf Jean Bertoin$^{(1)}$ and Alain Rouault$^{(2)}$}
\vskip 1cm
\noindent
\noindent
(1) {\sl Laboratoire de Probabilit\'es et Mod\` eles Al\'eatoires
and Institut universitaire de France,
 Universit\'e Pierre et Marie Curie,  175, rue du Chevaleret,
 F-75013 Paris, France.}
\vskip 2mm
\noindent
(2) {\sl {LAMA, B\^atiment Fermat,
Universit\'e de Versailles F-78035 Versailles. 
}}
\vskip 10mm

\[\mbox{{\bf Summary}}\] 
For a subcritical Galton-Watson process $(\zeta_n)$,  it is well known that under an $X \log X$ condition, 
the quotient $P(\zeta_n > 0)/ E\zeta_n$ has a finite positive limit. There is an analogous result 
for a (one-dimensional) supercritical branching random walk: when $a$ is in the so-called subcritical speed area,    
the probability of presence around $na$ in the $n$-th generation is asymptotically proportional to the corresponding expectation. 
In \cite{Rou93} this result was stated under a natural $X \log X$ assumption on the offspring point process and a
(unnatural) condition on the offspring mean. Here we 
prove that the result holds without this latter condition, in particular we allow an infinite mean and a dimension $d \geq 1$ for the state-space. 
As a consequence 
the result holds also for homogeneous fragmentations as defined in 
\cite{Bertoin01}, using the method of discrete-time skeletons; this completes the proof of Theorem 4 in \cite{BertRou2}. 
Finally, an application to conditioning on the presence allows to meet again the  probability 
tilting and the so-called additive martingale.
 \vskip 5mm
\noindent
{\bf A.M.S. Classification.}  {\tt 60 J 25, 60 F 10. }
\vskip 3mm
\noindent{\bf e-mail.} {\tt $(1):$ jbe@ccr.jussieu.fr , $(2):$
rouault@math.uvsq.fr }

\section{Introduction}
The common feature of many branching models consists in exponential growth. Let us recall some basic facts 
about a Galton-Watson process
$\zeta_n$  with finite mean started from $\zeta_0 = 1$. We have $`E(\zeta_n) = `E(\zeta_1)^n$ and

- if $`E(\zeta_1) > 1$ and $`E(\zeta_1 \log_+ \zeta_1) < \infty$, then
$\lim_n \zeta_n / `E(\zeta_n) = W > 0$ a.s. conditionally on non-extinction.

- if $`E(\zeta_1) < 1$ and $`E(\zeta_1 \log_+ \zeta_1) < \infty$, then,
\begin{equation}
\label{GW}
\lim_n \frac{`P(\zeta_n \geq 1)}{`E(\zeta_n)} = K > 0\,,
\end{equation}
 and there is a representation formula for the constant $K$ (see \cite{AsHer}). Moreover the condition $`E(\zeta_1 \log_+ \zeta_1) < \infty$ 
 is in some sense necessary.

In a discrete-time branching random walk (BRW), the initial ancestor is at the
origin  in
$\BBr ^d$ and the positions of its children form a point process $Z$. Each of  these
children has children in the same way: the positions of each family
relative to its parent is an independent copy of $Z$. Let $Z_n$ denote the point
process in
$\BBr ^d$ formed by the $n$-th generation. 

The intensity of $Z$ is the Radon
measure $\rho$ defined by 
\ben`E\left(\int_{\BBr ^d} f(x) Z(dx)\right)\,=\, \int_{\BBr ^d} f(x)
\rho(dx)\,.
\een 
and the intensity of $Z_n$ is the $n$-th convolution product $\rho^{*n}$.
We assume $1 < \rho(\BBr ^d) \leq +\infty$ (supercriticality) and set  
\ben
\label{zhat}\widehat Z (\theta) := \int_{\BBr ^d} \e^{\theta \cdot x} Z(dx)\quad
\hbox{and} \quad 
\Lambda (\theta) := \log `E (\widehat Z(\theta)) = \log \int_{\BBr ^d}
\e^{\theta\cdot x} \rho (dx)\,.
\een 
We shall also assume that $\Theta := \hbox{int}\ \{ \theta
: \Lambda(\theta) < \infty\} \not= \emptyset$.

The asymptotic behaviour as $n \rightarrow \infty$ of 
the 
%process
random measure $Z_n$, and more precisely, estimates
of $Z_n (na + I)$, where $a =\nabla\Lambda (\theta)$ and $I$ is some fixed
bounded set, have raised a considerable interest. 
This behaviour depends on  
 the
value of
$\theta$. 
The quantity 
$$ \Lambda^* (a) \,=\,\theta\cdot\nabla\Lambda(\theta) - \Lambda (\theta)$$
plays the role of the negative of the logarithm of the 
mean reproduction in the Galton-Watson process. 

In the range of supercriticality, i.e. for those  
$a \in \BBr ^d$ such that
 $ \Lambda^* (a) < 0$, 
the exponential rate of growth of $Z_n (na +I)$ is  
$- \Lambda^* (a)$ (see \cite{JDB77}) and a precise first order estimate  is given 
in \cite{JDB92} under an $X \log X$ type condition on $\widehat Z
(\theta)$.

In the range of subcriticality,  i.e. for those  
$a \in \BBr ^d$ such that
$ \Lambda^* (a) > 0$,  a precise analog of (\ref{GW}) i.e.
\ben
\lim_n \frac{`P(Z_n (na +I) \geq 1)}{`EZ_n (na +I)}= K(a) > 0
\een
 was proved in \cite{Rou93} under
certain conditions  which force in particular the finiteness of $\rho(\BBr^d)$,
and unfortunately this restricts the range of its applications. Our
 main purpose here is to show that the same estimate holds in fact under much weaker conditions.
Let us recall that the method is based on a representation formula for the presence probability by means of an auxiliary random walk. 
It is a discrete version of the Feynman-Kac formula used in  
probabilistic representations of solution of reaction-diffusion equations (\cite{BenRou93}, \cite{ChauRou88}, \cite{LaSe}).
Then a careful study of large deviations for functionals of this walk is needed, and our improvement takes place there.
  
Fragmentations, as defined in \cite{Bertoin01}, may be viewed as an extension of 
continuous-time branching random walks, but with a possible infinite offspring mean and infinite rate of branching.
 In \cite{BertRou2}, we succeeded to extend the results of the supercritical range to fragmentations, 
checking that the $X \log X$ assumption is automatically satisfied for the discrete-time skeleton. Moreover we proved that, in the subcritical range, the result 
can be deduced from its analogous version for BRW, by the method of discrete skeleton.

In Section \ref{brw} we state the result for BRW
%branching random walks 
 with a detailed proof in subsection \ref{proof}.
In Section \ref{frag} we extend the result to fragmentations, referring to Theorem 4 of \cite{BertRou2}. In subsection \ref{tilted} 
we show that conditioning on    $\{Z_t (at + I) \geq 1\}$ and letting $t \rightarrow \infty$ 
provides the tilted probability on the fragmentation process, already found in \cite{BertRou2}.
 
\section{Branching Random Walks}
\label{brw}
\subsection{Notation and representation formula}
To study spatial (Markov) branching processes, it is classical to use the probability generating functional.
We define its action on a Borel function $g$ satisfying $0\leq g\leq 1$ and such that $1-g$ has a compact support 
by the function 
$$F[g](x) = `E \exp \int_{\BBr^d} \log g(x+y) Z(dy)\ \,, \ x \in \BBr^d \,.$$
Its first moment is the linear operator acting on Borel functions with compact support by
$$M[g] (x) = `E \int g(x+y) Z(dy)\ \,, \ x \in \BBr^d \,,$$ 
(see \cite{AsHer} chap. V.1 and V.2).
Set $Z_n = \sum_i \delta_{z_i^n}$. Thanks to the Markov property, for every $n\geq 1$, the non linear operator $u_n$ defined by
\begin{equation}
\label{defu}
u_n[f] (x) = 1 - `E\left(\prod_i ( 1 - f(x + z_i^n))\right)\,,
\end{equation}
satisfies
$u_n[f] := 1 - F^n[1-f]$. Its linear counterpart, $v_n$ defined by
\begin{eqnarray}
\label{1}
v_n[f](x)= `E\left( \int_{\BBr ^d}  f(x+y) Z_n (dy) \right) = \int_{\BBr ^d} f(x+ z)
\rho^{*n} (dz)
\end{eqnarray}
satisfies $v_n[f] = M^n[f]$. We set $u_0[f]= v_0[f] = f$. 
When $f = {\bf 1}_J$ is the indicator function of a Borel set $J$, then
\begin{eqnarray}
u_n[{\bf 1}_J] (x) = `P (Z_n (J -x) \geq 1)\ \,, \ v_n[{\bf 1}_J] (x) = 
  `E (Z_n (J-x))\,.
\end{eqnarray}

%The evolution in mean is ruled by the convolution operator. For $n \geq 1$ and
%We start by introducing some notation.
%Let
% $f$  a  Borel function with compact support 
%Set  $v_0 =f$ and 
%For every $n \geq 1$, 
%we define the function $v_nf$ by:
%\begin{eqnarray}
%\label{1}
%v_nf (x) := `E\left( \int_{\BBr ^d} f(x+y) Z_n (dy)\right) \,=\, \int_{\BBr ^d} f(x+ z)
%\rho^{*n} (dz)\,.
%\end{eqnarray}
To study the propagation, let us introduce tilted (probability) measures\footnote{For simplicity we will 
assume that $\rho$ is not lattice and that its support is not contained in any lower-dimensional hyperplane.}. 
For $\theta \in \Theta$ and 
 $a = \nabla\Lambda(\theta)$, let $\big( S_n := \xi_1 + \cdots + \xi_n, n\in\BBn\big)$  be a
random walk with step 
distribution $\rho_\theta$ defined by
\begin{equation}
\label{tilt}
\int_{\BBr ^d} g(x) \rho_\theta (dx) = \int_{\BBr ^d} g(x- a)\e^{\theta \cdot x - \Lambda
(\theta)} \rho(dx)\,.
\end{equation}
The distribution $\rho_\theta$ is centered, 
its covariance matrix is Hess $ 
\Lambda(\theta)$ and we set $\sigma_\theta^2 = \det \hbox{Hess} \Lambda(\theta)$.
In the sequel, notation such as $`E_\theta$ will refer to expectations with respect to
this random walk.

From (\ref{1}) and (\ref{tilt}) we have the representation: 
\begin{eqnarray}
\label{repb2}
v_n[f] (x) =  \e^{-n \Lambda^*(a)} `E_{\theta} \left[\e^{- \theta \cdot S_n}f(x + S_n + na
)\right]\,.
\end{eqnarray}
Let us define the auxiliary function
%In the sequel, we use the notation
 $f_\theta$  by
\ben
f_\theta (y) = \e^{-\theta \cdot y}f(y)\,.
\een  
We have
\ben
\label{tlloc}
v_n[f](-na + c) = \e^{-n\Lambda^* (a)} \e^{\theta \cdot c}\,`E_\theta f_\theta (S_n + c)
\een
so that the local central limit theorem
entails that if $f_\theta$ 
%defined by $f_\theta (y) = \e^{-\theta \cdot y}f(y)$
 is directly Riemann integrable (DRI), then for every $c \in \BBr ^d$
\begin{equation}
\label{lclt0}
\lim_n \sigma_\theta (2\pi n)^{d/2}  \ \e^{n\Lambda^* (a)}v_n[f] (-an +c) 
= \e^{\theta \cdot c}\int_{\BBr ^d} 
f_\theta (y)
 dy\,. 
\end{equation}
%uniformly in  $x \in \BBr ^d$. %Up to a change of notation, the operator $u_n$ is 
To study the ratio between $u_n[f]$ 
%(which is a non-linear functional of $f$) 
and $v_n[f]$ 
%(which is a linear functional of $f$)
 we need  
the following representation formula (see \cite{Rou93}, formulas 4.2 and 4.3).

\begin{lemma}{\bf(\cite{Rou93})}
\label{lemrepres} Let $(P^{!x}, x \in \BBr ^d)$ be the family of reduced Palm distributions
defined by the disintegration formula:
\begin{equation}
\label{desint1}
`E\left(\int_{{\BBr ^d}} F(x , Z - \delta_x) Z(dx)\right)\,=\, \int_{\BBr ^d}\rho (dx) 
E^{!x}(F(x, Z))\,.
\end{equation}
and let for $r \geq 1$, $y, s \in \BBr^d$
\begin{equation}
\label{defh}
H_r [f](y, s) := E^{! y}\left( \int_0 ^1 \exp \langle Z, \log 
\left( 1 - \beta u_{r-1}[f] (s+.) \right) \rangle \, d\beta\right)\,,
\end{equation}
Then
\begin{eqnarray*}
\label{repres3}
u_n[f](x) = \e^{-n \Lambda^* (a)} `E_\theta\left\{\e^{-\theta \cdot
S_n} f(x +
S_n + na)
\prod_{1 \leq r\leq n}  
H_{r}[f]\left(a +  \xi_r, 
x +S_n - S_r +(n-r)a 
\right)
\right\}\,.  
\end{eqnarray*}
\end{lemma}
%It is a discrete version of the Feynman-Kac formula used in  
%probabilistic representations of solution of reaction-diffusion equations (\cite{BenRou93}, \cite{ChauRou88}).
We are now able to state the main result of this section.

\subsection{Main result}
\label{proof}
\begin{theorem}
\label{presdiscr}
Let $\theta \in \Theta$ satisfy
$\theta \cdot \Lambda'(\theta) - \Lambda (\theta) > 0$ 
 and 
$`E\left(\widehat Z(\theta)\log^{1+\varepsilon}\left(1+\widehat
Z(\theta)\right)\right)<\infty$ for some $\varepsilon>0$, and set 
  $a = \nabla\Lambda(\theta)$.  
Let  $f$ be a function with compact support, Riemann integrable, satisfying $0\leq f \leq 1$ and $\int_{\BBr^d} f(y) dy > 0$.
% with compact support.
 Then for every $c \in \BBr ^d$
\begin{equation}
\label{usurv}
\lim_{n \rightarrow \infty}\frac{u_n[f] (-an+c)}{v_n[f] (-an +c)} = 
\frac{\int_{{\BBr ^d}}  \e^{-\theta \cdot y}f(y)\,
 G[f](\theta, y)\, dy}{\int_{{\BBr ^d}}
\e^{-\theta \cdot y}f(y)
\, dy} 
\end{equation}
where
$$ G[f](\theta, y) \,=\, `E_\theta \left\{\prod_{r <\infty} H_r[f] (a + \xi_r, y-ar -
S_r)\right\}\,>\,0\,.$$
\end{theorem}

(One can see similar formulas for a branching diffusion in 
\cite{ChauRou88} Th. 1, \cite{LaSe} Prop. 1), \cite{BenRou93} Th.3. 3, \cite{Rouq} section 2.3). 
\medskip

\proof
In \cite{Rou93}, the proof was in two parts. The  first part 
was devoted to establish a lemma (4.1 p.33) 
and 
the second part consisted in checking conditions of the lemma. 
For easier reading we give here a complete proof, taking the main ideas of \cite{Rou93}, 
although the essential modifications 
take place in the part ``end of the proof'' below.

Let us first remark that if we set $\tau_c f := f(c+ \cdot)$ then $\tau_c$ commutes with $u_n$ and $v_n$. Moreover
$H_r[\tau_cf](y,s) = H_r[f](y, s+c)$, so that $G[\tau_cf](\theta , y) = G[f](\theta, y +c)$ and it is enough to prove (\ref{usurv}) for $c=0$.

Since $f$ is bounded by and has a compact support, there is 
some $b \in \BBr $ such that $f(x) = 0$ for $\theta \cdot x < b$.
Since $f \leq 1$
%If we set $f_\theta (x)  = \e^{-\theta \cdot x}f(x)$
we have 
\begin{equation}
\label{ftheta}
f_\theta(x) \leq \e^{-b} {\bf 1}_{\theta \cdot x \geq  b}\,.
\end{equation}
From Lemma \ref{lemrepres} we deduce 
\begin{eqnarray}
u_n[f] (-na) = \e^{-n \Lambda^* (a)} \BBe_\theta \left\{ f_\theta(
S_n )
\prod_{1 \leq r\leq n}  
H_{r}\left(a +  \xi_r, 
S_n - S_r -ra 
\right)
\right\}\,.  
\end{eqnarray}
Let us define
\begin{eqnarray*}
I_n^j &:=& \sigma_\theta (2\pi n)^{d/2} \BBe_\theta \big\{ f_\theta (S_n) \prod_{r=1}^j H_r
( a + \xi_r , S_n - S_r - ar) \big\} \ \ \hbox{for} \ \ 1 \leq j \leq n\\ I_\infty  ^j &:=&
\int_{\BBr ^d} f_\theta (s) \BBe_\theta \big\{  \prod_{r=1}^j H_r ( a + \xi_r , s - S_r - ar) 
\big\} \ ds \ \ \hbox{for} \ \ 1 \leq j\,.
\end{eqnarray*}
Taking into account (\ref{repb2}) and (\ref{lclt0}), the proof of (\ref{usurv}) 
can be reduced to showing that
\begin{equation}
\label{aim}
\lim_n I_n ^n = I_\infty ^\infty\,.
\end{equation}
A conditioning on $\xi_{j+1}, .. , \xi_n$ gives $I_n^j = \sigma_\theta (2\pi n)^{d/2}
\BBe_\theta [g_j (S_{n-j})]$ where
$$g_j (s) := \BBe_\theta \big\{ f_\theta (s + S_j) \prod_{k=1}^j H_k (a + \xi_k , s + S_j
-S_k -ka)\big\}\,.$$ We check that $g_j$ is DRI, so that the 
 local central
 limit theorem
yields
\begin{equation}
\label{lclt1}
\lim_n 
 I_n ^j = \int_{\BBr ^d} g_j (s) ds = I_\infty ^j\,,
\end{equation}
where the last equality comes from Fubini's theorem.
 Since $I_\infty ^j$ decreases to $I_\infty^\infty$ as $j\to\infty$, it is enough to prove
\begin{equation}
\label{aprouver}
\lim_j \limsup_{n > j} (I_n ^j - I_n ^n) = 0\,.
\end{equation}
Let us remark that we have no control on the uniformity in $j$ of (\ref{lclt1}).
If we search for a convenient upperbound for $I_n ^j - I_n ^n$, 
the difficulty comes from the term $S_r$ in $H_r (a + \xi_r , S_n - S_r -ar)$. Actually we
will give a lowerbound of $H_r$ after restricting the space. 

Let $\eta \in (0, \Lambda^* (a))$. 
For $i < k$, let $A_k^i := \{\theta \cdot S_\ell \leq -\ell\eta \ \hbox{for some} \ \ell \in (i+1 , k)\}$. 
Since $\theta \in \Theta$, there exists $\phi > 0$ such that $\int_{\BBr ^d} \e^{-\phi \theta\cdot x} \rho_\theta (dx) < \infty$, hence we can find $C_1$ and $C_2$ such that
for every $i$ and $k > i$
\begin{equation}
\label{LDP}
\BBp_\theta (A_k^i) \leq C_1 \e^{-C_2 i}\,.
\end{equation}
From now, $C_k , k \geq 1$ denote strictly positive finite constants depending on $\theta, a, b, \eta$ but not on other quantities.
\medskip

\noindent{\bf First step: Restriction of the space}
\begin{eqnarray}
\label{entrois}
\nonumber
0 < \frac{I_n ^j - I_n ^n}{\sigma_\theta(2\pi n)^{d/2}} &\leq&  \BBe_\theta\big[f_\theta
(S_n)\big(1- \prod_j^n\big)\big]\leq \\  &\leq& \BBe_\theta\big[f_\theta (S_n)
; A_n ^
j\big]
+ \BBe_\theta\big[f_\theta (S_n)\big(1- \prod_j^{n}\big); (A_n ^j)^c \big] 
\end{eqnarray}
Introducing an intermediate integer $n_1 \in (j, n)$ we have (see (\ref{ftheta}))
\begin{eqnarray}
\label{truc}
\BBe_\theta\big[f_\theta (S_n)
; A_n ^
j\big] \leq 
e^{-b} \BBp_\theta(A_n ^{n_1} ) 
+\BBe_\theta\big[f_\theta (S_n)
; A_{n_1}^j 
\big]
\end{eqnarray}
and, conditioning on 
$S_k , k \leq n_1$,
\begin{eqnarray*}
\BBe_\theta\big[f_\theta (S_n)
 ; A_{n_1}^{j} \big] 
= \BBe_\theta\big[\int_{\BBr ^d} f_\theta (S_{n_1} +s)\rho_\theta^{* (n-n_1)} (ds)  ; A_{n_1}^{j}\big]\,.
\end{eqnarray*}
Setting, for $q \geq 1$
\begin{eqnarray}
\label{granddelta}
\Delta_q := \sup_s \big|\sigma_\theta(2\pi q)^{d/2} \BBe_\theta f_\theta (s + S_q) - \int_{\BBr ^d} f_\theta (x) \exp- \frac{\Vert x-s\Vert^2}{2q\sigma_\theta^2} dx\big|\, 
\end{eqnarray}
we get
\begin{eqnarray*}
\sigma_\theta [2\pi (n-n_1)]^{d/2}\int_{\BBr ^d} f_\theta (S_{n_1} +s)\rho_\theta^{* (n-n_1)} (ds) \leq \Delta_{n-n_1} + \int_{\BBr ^d} f_\theta (x) dx
\end{eqnarray*}
hence
\begin{equation}
\label{trucbis}
\sigma_\theta [2\pi n]^{d/2}\BBe_\theta\big[f_\theta (S_n)
 ; A_{n_1}^{j} \big]\leq \big(1 - \frac{n_1}{n}\big)^{-d/2} \big\{\Delta_{n-n_1} + \int_{\BBr ^d} f_\theta (x) dx \big\} P_\theta (A_{n_1}^{j})\,.
\end{equation}
Fix $j$. Choosing $n_1$ such that $n^{d/2} \e^{-C_5 n_1} \rightarrow 0$, we get from (\ref{truc}),(\ref{LDP}), (\ref{trucbis})  and the local central limit theorem (\cite{St} Lemma 2) :
\begin{equation}
\label{premier}
\limsup_{n > j} \sigma_\theta (2\pi n)^{d/2} \BBe_\theta\big[f_\theta (S_n)
; A_n ^
j\big] \leq C_1 \e^{-C_2 j}\,.
\end{equation}
\medskip

\noindent{\bf Second step : lowerbound for $H_r$}

We start with 
\begin{equation}
\label{u<v}
u_n[f](x) \leq \min \left( v_n[f](x) , 1\right)\ \,, \ x \in \BBr^d \,, \ n\geq 0\,.
\end{equation}
%From (\ref{ftheta}) and (\ref{repb2}) we have
Since $f(x) \leq \e^{\theta\cdot x -b}$
we have
%and the exponential Chernov bound (see (\ref{ftheta}) and (\ref{repb2}))
\begin{equation}
\label{chernov}
v_n[f](x) \leq \e^{n\Lambda(\theta)} \e^{\theta \cdot x - b}\ \,, \ x \in \BBr^d \,, \ n\geq 0\,.
\end{equation}
Fix $y$ satisfying $\theta \cdot y \geq  - r\eta$. Applying (\ref{chernov}) 
with  $x = -y - ar +z +z_i$ 
we see that 
\begin{equation}
\label{ineqv} 
v_{r-1}[f] (-y -ar + z + z_i)\leq \e^{-r [\Lambda^* (a) - \eta]  + \theta \cdot (z + z_i)  -b -\Lambda (\theta)}
\end{equation}
for every $r \geq 1$ and $i\geq 1$.
Since $\Lambda (\theta) < \infty$, the random variable $\widehat Z (\theta)$ defined by (\ref{zhat}) is 
$`P$-a.s. finite, 
hence
\begin{equation}
\label{inegtau}
\log \widehat Z (\theta) %= \log \int_{\BBr ^d} \e^{\theta \cdot x} Z(dx) 
\geq 
\sup_i \{\theta \cdot z_i \} =: \tau(Z)  \,.
\end{equation}
Let 
$B_r := \{ Z : \tau(Z) < r[\Lambda^* (a) - \eta] +\Lambda (\theta) - \theta  \cdot z  +
b\}$. 
For  $Z \in B_r$,  the bound in (\ref{ineqv}) is less than or equal
to $1$  so that (\ref{u<v}) gives 
\begin{eqnarray*}
u_{r-1}[f](-y -ar + z + z_i) \leq 
\e^{-r [\Lambda^* (a) - \eta]  + \theta \cdot (z + z_i)  -b -\Lambda (\theta)}
\leq 1
\end{eqnarray*}
and then 
(when $\theta \cdot y \geq  - r\eta$)
\begin{eqnarray}
\label{16}
H_r (x , -y -ar +z)\geq {\cal H}_r (x, z) := 
E^{!x}\left({\bf 1}_{B_r}\int_0 ^1 \e^{ 
- w(r, Z,\beta, z)} d\beta\right)
 \end{eqnarray}
where
\begin{eqnarray}
\label{17}
w(r, Z,\beta , z) := - 
\int_{\BBr ^d} \log \left( 1 - \beta 
\e^{-r [\Lambda^* (a) - \eta]  + \theta \cdot (z + \zeta)  -b -\Lambda (\theta)}
\right)\, Z(d\zeta)
\,.
\end{eqnarray}

\noindent{\bf End of the proof}

Coming back to the second term of (\ref{entrois}) and using (\ref{16}) we get 
\begin{eqnarray}
\label{tructer}
\nonumber
 \BBe_\theta\Big[f_\theta (S_n)
\big(1- \prod_{j}^n\big)\!&;&\!\big(A_n^{j})^c \Big] \leq  \BBe_\theta \Big[f_\theta (S_n) 
\sum_{k=j}^n \big(1 -{\cal H}_k (a+ \xi_k , S_n)\big)\Big]\\
&=& \BBe_\theta \int_{\BBr ^d} f_\theta ( S_{n-1} +x ) \big[\sum_{k=j}^n \big(
1 - {\cal H}_k (a + x , S_{n-1} + x  )\big) \big] \rho_\theta (dx)
\end{eqnarray}
We want to bound $\sum_{r=j}^n \big(1 - {\cal H}_r (x,z)\big)$.
Adding up, from (\ref{16}) and (\ref{17})
\begin{eqnarray}
\label{1+2}
\nonumber
\sum_j ^n [1 -{\cal H}_r(x, z)] &\leq& E^{!x}\left( \sum_j ^n {\bf 1}_{B_r ^c} (Z)
\right) + E^{!x}\left( \int_0 ^1 \big(\sum_j ^n {\bf 1}_{B_r}\big(1 -\e^{- w(r,Z, \beta,
z)}\big) \big) \, d\beta \right)\\ &:= & J_1 (x,z)+ J_2 (x,z)\,.
\end{eqnarray}
From the definition of $B_r$ we have
\begin{eqnarray*}
\sum_j ^n {\bf 1}_{B_r ^c} (Z) \leq \left(\frac{ \tau(Z) + \theta \cdot z -b - \Lambda(\theta)}{\Lambda^* (a) - \eta} -j +1\right)^+\,,
\end{eqnarray*}
which, from inequality (\ref{inegtau}), gives by integration
\begin{equation}
\label{boundi1}
J_1 (x, z)\leq 
E^{!x}\left( \left(\frac{ \log \widehat Z(\theta) + 
\theta \cdot z -b - \Lambda(\theta)}{\Lambda^* (a) - \eta} -j +1\right)^+\,
\right)
\end{equation}
It remains to give an upperbound for $J_2 (x, z)$.
The function $r \mapsto w(r, Z, \beta,z) 
$
is decreasing in $r$, so by the classical sum-integral comparison, we obtain, for $Z \in B_j$
\begin{eqnarray*}
\sum_j ^n [1 - \e^{-w(r,Z,\beta, z)}] \leq \int_s ^\infty [1 - \e^{-w(r,Z,\beta,z)}] dr \leq C_3 \log [1 + w(s, Z, \beta, z)]
\end{eqnarray*}
where the last inequality comes from
$-\frac{\partial w}{\partial r} \geq w [\Lambda^* (a)-\eta] \,.$

The inequality\footnote{It can be proved using concavity of logarithm and Jensen's inequality.}:
\begin{eqnarray*}
\int_0^1 \log \left[1 - \sum_i \log (1 -\beta a_i)\right] d\beta \leq  \log \big(1 + \sum_i a_i \big)
\end{eqnarray*}
for $0\leq a_i \leq 1, i=1, \ldots$ gives
\begin{eqnarray}
\label{boundi2}
J_2 (x ,z)
\leq C_4 E^{!x}\left( \log ( 1 + \widehat Z(\theta) \e^{
- j [\Lambda^* (a) -\eta]  + \theta \cdot z -b -\Lambda(\theta)
})\right)\,.
\end{eqnarray}
Combining (\ref{boundi1}) and (\ref{boundi2}) gives
\begin{eqnarray}
\label{above}
J_1 (x ,z)+ J_2 (x ,z) \leq C_5 E^{!x}\left( \log ( 1 + \widehat Z(\theta) \e^{
- j [\Lambda^* (a) -\eta] + \theta \cdot z -b  -\Lambda(\theta)
})\, \right)\,.
\end{eqnarray}
Setting
$$A(x, z) = E^{!x}\left( \big[ \log 1 + \widehat Z(\theta) \e^{\theta \cdot z
-b-\Lambda(\theta)}\big]^{1+\epsilon} \right)$$
and
applying inequality (A2) p.38 of \cite{Rou93} we see that for every $\epsilon > 0$, 
the right hand side of (\ref{above}) is bounded by $C_3 j^{-\epsilon} A(x, z)$ so that, from (\ref{1+2}),
\begin{eqnarray*}
\sum_j ^n[1 - {\cal H}_r (x, z)] \leq C_6 j^{-\epsilon} A(x, z)\,.
\end{eqnarray*}
Thanks to  (\ref{tructer}) that entails
\begin{eqnarray*}
 \BBe_\theta\big[f_\theta (S_n)
\big(1- \prod_{j}^n\big) ; \big(A_n^{j})^c \big] 
\leq C_6 j^{-\epsilon}\BBe_\theta \left(B(S_{n-1})\right)
\end{eqnarray*}
where
\begin{eqnarray*}
B(s) := \int_{\BBr ^d} f_\theta(x+s) A(a+x, s+x)\, \rho_\theta (dx)\,. 
\end{eqnarray*}
It should be clear that the function $B$ has bounded variation. 
By Fubini's theorem and the disintegration formula (\ref{desint1}), (recalling that $f$ has a compact support included in $\{y : \theta \cdot y \geq b\}$) we get 
\begin{eqnarray*}
\int_{\BBr ^d} B(s) ds \leq C_5 \ \BBe \left(\int_{\e^{  - \Lambda(\theta)}} ^\infty v^{-2} 
\big[ \log 1 + \widehat Z(\theta) v
\big]^{1+\epsilon} \widehat Z(\theta)  dv\right)\,. 
\end{eqnarray*}
Since $$\log (1 + \widehat Z(\theta) v) \leq \log (1 + v)  + \log_+ \widehat Z(\theta)$$
and by convexity
$$[\log (1 + \widehat Z(\theta) v)]^{1 +\epsilon} \leq C_7 \big([\log (1 + v)]^{1 +\epsilon} + 
[\log_+ \widehat Z(\theta)]^{1 +\epsilon}\big)\,,$$
we see that $B$ is integrable under the assumptions of Theorem \ref{presdiscr}.

Invoking again the local central limit theorem,
\begin{eqnarray*}
\lim_n \ \sigma_\theta (2\pi n)^{d/2} \BBe_\theta \left( B(S_{n-1})\right) = \int_{\BBr ^d}
B(s) ds\,,
\end{eqnarray*} so that from (\ref{tructer})
\begin{equation}
\label{dernier} \limsup_{n>j} \sigma_\theta (2\pi n)^{d/2} \BBe_\theta\big[f_\theta (S_n)
\big(1- \prod_{j}^n\big) ; \big(A_n^{j})^c \big] 
\leq C_3 j^{-\epsilon} \int_{\BBr ^d} B(s) ds\,.
\end{equation}
Taking into account (\ref{premier}) and (\ref{dernier}), we see that (\ref{aprouver}) holds, which entails (\ref{aim}) and 
completes the proof of  (\ref{usurv}). \QED 

\section{Fragmentations}
\label{frag}
\subsection{Notations and main result}

We follow the notations of \cite{Bertoin01}, \cite{Beres02}, \cite{Bertoin03} and \cite{BertRou2} and refer to these papers for details. We  work
with the space of numerical sequences
$$\s\,:=\,\left\{{\bf s}=(s_1,\ldots): s_1\geq s_2\geq\cdots \geq0\hbox{ and }
\sum_{1}^{\infty}s_i\leq1\right\}\,,$$
which should be thought as the set of ranked masses of the
fragments resulting from the split of some object with unit total mass. 
We consider a family of Feller processes $X=(X_t, t\geq0)$ with
values in $\s$ and c\` adl\` ag paths.
For every $a\in[0,1]$, we let $`P_a$ denote the law
of $X$ with initial distribution $(a,0,\ldots)$ (i.e. the process
starts from a single fragment with mass $a$). We say that $X$ is a
 (ranked) {\it homogeneous fragmentation} if the following two
properties hold: 

\noindent $\bullet$ (Homogeneity property) For every $a\in[0,1]$, the law
of 
$aX$ under $`P_1$ is $`P_a$.

\noindent $\bullet$ (Fragmentation property) For every ${\bf s}=(s_1,
\ldots)\in\s$, the process started from $X(0)={\bf s}$ can be obtained as
follows. Consider $X^{(1)},
\ldots$ a sequence of independent processes with respective laws
$`P_{s_1}, \ldots$, and for every $t\geq0$, let $\hat X(t)$ be the random
sequence obtained by ranking in decreasing order the terms of the
random sequences $X^{(1)}(t), \ldots$. Then $\hat X$ has the law of $X$
started from ${\bf s}$. Here we consider homogeneous fragmentations with no erosion (\cite{Bertoin03}), so that 
the distribution of $X$ 
% consider The distribution of the process $X$ 
 can be characterized  by a 
measure $\nu$ on $\s$, called the {\it dislocation measure}. Informally
$\nu$ specifies the rates at which a unit mass splits. It has to fulfil the conditions
$\nu(\{(1,0,\ldots)\})=0$ and 
\begin{equation}
\label{1-s1}
\int_{\s}\left(1-s_1\right)\nu(d{\bf s})\,<\, \infty\,.
\end{equation}
 We shall assume  that
\begin{equation}\label{pasperte}
\nu\left(\left\{{\bf s}\in\s:
\sum_{i=1}^{\infty}s_i<1\right\}\right)\,=\,0\,,
\end{equation}
 which means that no
mass is lost when a sudden dislocation occurs, and more precisely,
entails that the total mass is a conserved quantity for the fragmentation
process (i.e. $\sum_{i=1}^{\infty}X_i(t)=1$ for all $t\geq0$, 
$`P_1$-a.s.). Moreover we also  exclude the
trivial case when $\nu\equiv 0$.  
% In other words, the fragmentation does not produce dust.
%In the sequel, we shall also implicitly exclude the
%trivial case when $\nu\equiv 0$.  

Given a real number $r>0$, we say that a dislocation measure $\nu$ is
$r$-geometric  if it is finite and is carried by the 
subspace of configurations ${\bf s}=(s_1,\ldots)\in\s$ such that
$s_i\in \{r^{-n}, n\in\BBn\}$. This holds if and only if
$`P_1(X_i(t)\in\{r^{-n}, n\in\BBn\}
\hbox{ for every }i\in\BBn)=1$ for all $t\geq0$.
We say that $\nu$ 
is {\it non-geometric} if it is not $r$-geometric for any $r>0$.

The empirical measure of the logarithms of the fragments
\ben
Z^{(t)} (dy) := \sum_{i=1}^\infty \delta_{\log X_i (t)}, \ t\geq 0
\een
can be viewed as the generalization  of a branching process in continuous time, 
with a possible infinite offspring mean and infinite rate of branching, corresponding to $\nu(\s) = \infty$.
If we define
\ben
\Phi(p) := \int_{\s} \left(1 - \sum_i x_i^{p+1} \right) \nu(d{\bf x}) , \ p > \underline p
\een
where
$$\underline p\,:=\,\inf\left\{p\in\BBr:
\int_{\s}\sum_{i=2}^{\infty}s_i^{p+1}\nu(d{\bf s})<\infty\right\}
\,,$$
then 
\ben
`E\left(\int_{\BBr^d} e^{(p+1)y}Z^{(t)} (dy)\right) = `E\left(\sum_i X_i (t)^{p+1}\right) = \exp (-t \Phi(p))
\een
% the counting process and, if $\Phi$ is the Laplace exponent of the corresponding subordinator, let
and \begin{eqnarray}
M(p, t) := \int_\BBr e^{(p+1)y + t \Phi(p)} Z^{(t)} (dy) = e^{t\Phi(p)} \sum_i X_i^{p+1}(t)\,,  \ \ t\geq 0 
\end{eqnarray}
is the so-called additive martingale. The function $\Phi$ is concave, analytic and increasing. 
It is the Laplace exponent of a subordinator ($\chi_t$)
\ben
\label{subord}\exp (-t \Phi(p)) = `E \exp (-p \chi_t)
\een (see \cite{Bertoin01} for details).
If $\bar p$ denotes the unique solution of the equation
$$\Phi(q)\,=\,(q+1)\Phi'(q)\,,\qquad q>\underline p\,,$$
we have 
\begin{eqnarray*}
\Phi(q) -(q+1)\Phi'(q) &<& 0\ \ \hbox{for}  \ \underline p < q < \bar p\\  
\Phi(q) -(q+1)\Phi'(q) &>& 0\ \ \hbox{for}  \  p > \bar p\,.
\end{eqnarray*}
Since $\log X_1(t)$ (logarithm of the maximal size) grows as $t\rightarrow \infty$ like $-t \Phi'(\bar p)$, 
we say that $\{a = -\Phi'(p) ; \ \underline p < p < \bar p\}$ is the supercritical range 
and $\{a = -\Phi'(p) ; \ \bar p < p \}$ the subcritical range. 

Let us fix $\alpha < \beta$ and $p > \bar p$. Here we are interested in the asymptotic behaviour of
\begin{eqnarray}
\label{defUV}
U(t, x) &:=& `P( Z^{(t)}\big([\alpha + x, \beta +x]\big) \geq 1) 
\\
\nonumber
V(t, x) &:=& `E Z^{(t)} \big([\alpha +x ,  \beta +x]\big)\,, 
\end{eqnarray}
for $x = -t \Phi'(p)$ and $t\rightarrow \infty$.
%Let $f = \BBone_{[\alpha , \beta]}$ and $\theta = p+1$.
%We recall that the dislocation measure $\nu$ of the fragmentation satisfies
%\begin{equation}
%\label{1-s1}
%\int_{\cal S} (1 - x_1) \nu (d{\bf x}) < \infty\,.
%\end{equation}

The following theorem is proved in \cite{BertRou2} using time-discretization and taking for granted Theorem 2 of the present paper. 
\begin{theorem}{\bf(\cite{BertRou2})}
 Assume that the dislocation measure $\nu$ is
non-geometric. 

\label{presfrag}
\begin{itemize}
\item[{\rm (i)}] If $p > \underbar p$, we have 
\begin{eqnarray*}
\lim_{t\to\infty} 
\sqrt t \, \e^{-t((p+1)\Phi'(p)-\Phi(p))} 
V(t, -t\Phi'(p))  \,=\, \frac{1}{\sqrt{2\pi
|\Phi''(p)|}} 
(p+1)^{-1}
 \left(\e^{-(p+1) \alpha}-\e^{-(p+1) \beta}\right)
\,. 
\end{eqnarray*}
\item[{\rm (ii)}] If $p > \bar p$, there exists a positive finite constant
$K_{p}$ such that
\begin{eqnarray*}
\label{usurvfg}
\lim_{t \rightarrow \infty}\frac{U(t, -t\Phi'(p))}{V(t, -t\Phi'(p) )} = K_p \,. 
\end{eqnarray*}
\end{itemize}
\end{theorem}
\medskip

\noindent{\bf Remark}
Actually, for fixed $c \in \BBr$, we have 
\begin{eqnarray}
\label{fixedc}
\lim_{t\to\infty} 
\frac{V(t, -t\Phi'(p) +c)}{V(t, -t\Phi'(p))}  \,=\, \e^{(p+1)c}
\ \ , \ \  
%\end{eqnarray*}
%and
%\begin{eqnarray*}
\lim_{t\to\infty}\frac{U(t, -t\Phi'(p)+c)}{V(t, -t\Phi'(p)+c)} = K_p 
\end{eqnarray}
(exactly as in (\ref{usurv})).

\subsection{Comment}
 In \cite{BertRou2} it was shown that $K_p$ may be obtained as the constant coming from any 
discrete skeleton (i.e. for any arbitrary choice of the time mesh). Actually, a careful analysis may give a representation formula for the constants in 
terms of an underlying L\'evy process instead of a (skeleton) random walk. We state it without proof not to overburden the paper. We have 
\begin{equation}
\label{usurvf}
\lim_{t \rightarrow \infty}\frac{U(t, -at)}{V(t, -at)} = 
\frac{\int_\alpha^\beta e^{-(p+1)y} G(p, y) dy}{\int_\alpha^\beta e^{-(p+1)y} dy} 
\end{equation}
where $G(p,y)$ is an expression that we explain now.

In the spirit of Section 2, we define
$$U_t[f](x)= 1 - `E\left(\prod_i (1 - f(x + \log X_i (t)))\right)\ \ , \ \ V_t[f](x)= `E \left(\sum_i f(x + \log X_i (t) )\right)$$
so that $U(t,x) = U_t[{\bf 1}_{[\alpha , \beta]}] (-x)$. 
The random walk $( S_n , n\geq1)$ under the law $`P_\theta$ in Section \ref{brw} is now replaced by a 
certain L\'evy process $(\hat\zeta_t)$. 

Recall that $\Phi$ is the Laplace exponent of the L\'evy process $(\chi_t)$ (see (\ref{subord})). Let $L$ be its characteristic measure  so that
$$\Phi(\lambda) = \int_{]0, \infty[} (1 - e^{-\lambda x}) L(dx)\,.$$ 
Then $(\hat\zeta_t, t \geq 0)$ is the dual of the L\'evy process
whose characteristic measure is
$\tilde L_{p} (d\Delta) = \e^{-p\Delta} L(d\Delta)$ and %
drift coefficient
%mean expectation
  $a$, so that its mean expectation is $0$ (see \cite{Bertoinb96}). We have in particular
\ben
\label{Levyconv}  
V_t[f](-at +c) = \e^{(t(p+1)\Phi'(p)-\Phi(p))} \e^{(p+1)c} \ `E f_{p+1}(\hat\zeta_t + c)\,.
\een
%It plays the role of the random walk $( S_n , n\geq1)$ under the law $`P_\theta$ in Section \ref{brw}. 

To study $U_t$ we need more definitions.
Let ${\bf x}_*$ be a `size biased pick' from the sequence ${\bf
x}$, i.e. a random variable with values in $]0,1[$ such that for every ${\bf x}\in \s$ and $i \in \BBn$
$$`P({\bf x}_* =x_i \,|\, {\bf x}= (x_1, \cdots)) = x_i$$
or in other words $`E(g({\bf x}_*)| {\bf x}) = \langle {\bf
x},\overline g\rangle$,
%$\BBe(g({\bf x}_*))=\langle {\bf
%x},\overline g\rangle$, 
where $\overline
g(y)=yg(y)$. At last we denote by 
 ${\bf x}^{!}$ the random sequence obtained by removing the 
size biased pick ${\bf x}_*$ from
%  the latter in
the sequence ${\bf x}$.

Let $m$ be 
the "distribution" under $\nu$ of the size biased pick,
i.e.
$$\int_{(0,1)}g(x)m(dx)\,=\,\int_{\s}\nu(d{\bf x})\BBe(g({\bf
x}_*)\,|\,{\bf x})\,.$$ We check easily that
$\int_{(0,1)}^{}(1-x)m(dx)<\infty$.
For
 $x\in(0,1)$, we denote by
$\nu^{!x}$
the probability measure defined on $\s$ by the
disintegration
\begin{equation}
\label{desint2}
\int_{\s}\nu(d{\bf x})\BBe\left(H({\bf x}_*,{\bf x}^!)\right)
\,=\,\int_{]0,1[}m(dx)\int_{\s}\nu^{!x}(d{\bf x})H(x,{\bf x})\,.
\end{equation}
In words,
 $\nu^{!x}$ is the conditional law under 
$\nu$ 
of ${\bf x}^!$
knowing ${\bf x}_*=x$. This definition is similar to (\ref{desint1}) in the BRW framework.

With that notation we may define
\begin{equation}
\label{defh1}
H_r ( \xi, z)
= \int_0 ^1\int_\s \exp 
\left(\sum_{i=1}^\infty \log (1 - \beta U_r(z + \log x_i))\right)
 \nu^{!
e^\xi
} (d{\bf x}) d\beta\,,
\end{equation}
and the function $G(p,y)$ in (\ref{usurvf}) is 
$$ G(p, y) = \BBe \left\{\prod_{r <\infty, \ \Delta\hat\zeta_r \not= 0} H_r (\Delta\hat\zeta_r, y-ar-\hat\zeta_r)\right\}\,.$$
where $\Delta\hat\zeta_r = \hat\zeta_r-\hat\zeta_r^-$.

\subsection{The tilted probability in the subcritical range}
\label{tilted}
This section is an easy extension to fragmentations of the corresponding result for branching Brownian motion \cite{ChauRou88} section 6, or branching random walk 
\cite{Rou93} p. 30-31 but or the sake of completeness we give the complete proof. For fixed $a$ in the subcritical range we show that 
conditioning on the presence of fragments of size $\e^{-at + 1}$ at time $t$
gives in the limit $t \rightarrow \infty$ the tilted probability $`P^{(p)}$ 
for the measure valued process $(Z_s , s\geq 0)$. 
This tilted probability, described in \cite{BertRou2} Section 3.3, 
 is the h-transform of the probability by the martingale $M(p, \cdot)$. This result is in the spirit of many results on conditioning 
 conditioning branching spatial processes and superprocesses and related to the notion of immortal particle due to Evans (see an extensive bibliography in \cite{EnKyp}).
  
 Let $({\cal F}(s) , s \geq 0)$ be the natural filtration defined by 
${\cal F}(s) := \sigma (Z_r , r\leq s)$ and ${\cal F}= \wedge {\cal F}_s$.
\begin{proposition} Let us fix $\alpha < \beta$, 
%$J = [\alpha, \beta]$.
  $p > \bar p$ and $a = -\Phi' (p)$. Then, for every $A \in {\cal F}$ we have:
%$a,p$ as in Theorem \ref{presfrag} b), 
$$\lim_{t \rightarrow \infty} `P(A \,|\  
Z_t ( [at + \alpha, at +\beta)
 \geq 1) =  {`P}^{(p)}(A) $$
where, for every $s > 0$,
$${`P}^{(p)} |_{{\cal F}(s)} = M(p,s) `P |_{{\cal F}_s}\,.$$
\end{proposition}
\medskip

\proof
Set $J = [\alpha, \beta]$. Fix $s > 0$ and  $G_s \in {\cal F}_s$. We have for $t > 0$
\begin{eqnarray}
\label{37}
`P\big(G_s \,|\ Z_{t+s} ( a(t+ s) + J) \geq 1)
= `E\left[G_s \frac{`P [Z_{t+s} (a(t+ s)+ J  ) \geq 1] \,|\ {\cal F}_s]}{
`P [Z_{t+s} ( a(t+ s) + J) \geq 1]
}\right]
\end{eqnarray}
From (\ref{defUV})
%With the notations of the above subsection, we have
\begin{eqnarray*}
`P [Z_{t+s} ( a(t+ s) + J) \geq 1] 
 = U_{t+s} ({-a(t+ s)}) 
\end{eqnarray*}
and, by the fragmentation property
\begin{eqnarray}
\label{num}
`P [Z_{t+s} (a(t+ s) + J) \geq 1 
 \,|\,{\cal F}_s] = 1 - \prod_j (1 - A_j (t))
\end{eqnarray}
where $A_j (t) := 
 U_{t}(\log X_j(s) -a(t+s))\,.$ 

To apply the results of the above section, we set  
$r_t := \sigma_p \sqrt{2\pi t}\ e^{t\Lambda^* (a)}$
(which tends to infinity with $t$)
and  $K'_p :=  \int_\alpha^\beta e^{-(p+1)y}  dy\,.$
 From Theorem \ref{presfrag} ii)
\begin{equation}
\lim_t r_t U_{t+s} (-a(t+ s)) =  
K_p K'_p \ e^{-s\Lambda^* (a)}\,.
\end{equation}
To handle the expression (\ref{num}), we apply Theorem \ref{presfrag} i) 
 for every $j$:
\begin{eqnarray*}
\lim_t 
r_t
 A_j (t)
= K_p K'_p\  X_j (s)^{p+1} e^{-(p+1) as} 
\,,
\end{eqnarray*}
but we need a uniform bound. 
%First, we can bound $U$ by $V$ from (\ref{Levyconv}) and 
From (\ref{Levyconv}) and the local central limit theorem, there exists $\epsilon_t \rightarrow 0$ such that, for every $c \in \BBr$
$$r_t V_t (-at + c) \leq K'_p \big( 1+ \epsilon_t\big) \e^{(p+1)c}\,.$$
This yields:
%and apply 
%again Theorem \ref{presfrag} i)  to get
\begin{eqnarray*}
r_t
 A_j (t)
\leq 
r_t
 V_{t}(\log X_j(s) -a(t+s))
\leq K'_p \big( 1+ \epsilon_t\big)  X_j (s)^{p+1} e^{-(p+1) as}\,.
\end{eqnarray*}
%where $\epsilon_t$ tends to zero and is independent of $X_j(s)$. 
By dominated convergence, we have a.s. 
$$\lim_t r_t \sum_j A_j (t) =  \ \sum_j X_j (s)^{p+1} e^{-(p+1) as}
 = K_p K'_p\  M(s, p) e^{-s \Lambda^* (a)}$$
 and since
$$\sum_j A_j (t) - \big(\sum_j A_j (t)\big)^2 \leq 1 - \prod_j (1 - A_j (t)) \leq \sum_j A_j (t)\,,$$
we get
\begin{eqnarray*}
\lim_t \frac{`P [Z_{t+s} ( a(t+s) + J) \geq 1 \,|\, {\cal F}_s]}{
`P [Z_{t+s} ( a(t+s) + J)   \geq 1]} = \lim_t \frac{r_t \sum_j A_j (t)}
{r_t U_{t+s} (-a(t+ s))}= M(s,p)
\end{eqnarray*}
a.s.. Invoking again the dominated convergence theorem allows to conclude from (\ref{37})
$$\lim_t  `P\big(G_s \ | \ Z_{t+s} ( a(t+ s) + J) \geq 1)
= `E [G_s M(p,s)]$$
which ends the proof of the proposition. \QED
\bibliographystyle{plain} 
\small
\bibliography{presence}
\end{document}